\documentclass[12pt]{amsart}
\usepackage{amssymb}
\usepackage{amsmath}
\usepackage{pictex}
\newcommand\0{\mathbf{0}}

\newcommand\AFF{\operatorname{\sf Aff}}
\newcommand\al{\alpha}
\newcommand\AND{\quad\mbox{and}\quad}
\newcommand\BB{\mathcal B}
\newcommand\bd{\partial}

\newcommand\CC{\mathcal C}
\newcommand\cf{\curlywedge}

\newcommand\de{\delta}
\newcommand\DL{\mathsf{DL}}

\newcommand\ep{\varepsilon}

\newcommand\fl{\mathfrak{f}}

\newcommand\Ga{\Gamma}
\newcommand\geo[1]{\overline{#1}}

\newcommand\HH{{\mathcal H}}
\newcommand\hor{\mathfrak{h}}
\newcommand\Hor{\operatorname{\sf Hor}}

\newcommand\KK{\operatorname{\sf K}}
\newcommand\la{\lambda}

   \newcommand\MM{\mathcal M}
\newcommand\om{\omega}
\newcommand\Om{\Omega}
\newcommand\mm{\mathsf{m}}
\newcommand\N{\mathbb N}

\newcommand\Prb{\operatorname{\sf Pr}}
   \newcommand\R{\mathbb R}

\newcommand\scs{\scriptstyle}
\newcommand\sign{\operatorname{sign}}

\newcommand\SH{{\mathcal S}}
\newcommand\supp{\operatorname{\sf supp}}
\newcommand\T{\mathbb T}

\newcommand\U{{\mathcal U}}
\newcommand\up{\mathfrak{u}}
\newcommand\wh{\widehat}
\newcommand\wt{\widetilde}
\newcommand\Z{\mathbb Z}

\numberwithin{equation}{section}

\newtheoremstyle{mythm}
  {9pt}
  {9pt}
  {\itshape}
  {0pt}
  {\bfseries}
  {}
  { }
  {\thmnumber{(#2)}\thmname{ #1}\thmnote{ #3}}

\newtheoremstyle{mydef}
  {9pt}
  {9pt}
  {\normalfont}
  {0pt}
  {\bfseries}
  {}
  { }
  {\thmnumber{(#2)}\thmname{ #1}\thmnote{ #3}}

\theoremstyle{mythm}
\newtheorem{thm}[equation]{Theorem.}
\newtheorem{pro}[equation]{Proposition.}
\newtheorem{lem}[equation]{Lemma.}

\theoremstyle{mydef}

\newtheorem{rmk}[equation]{Remark.}

\begin{document}
\title{\large Positive harmonic functions for\\ 
semi-isotropic random walks on trees, \\
lamplighter groups, and DL-graphs} 
\author{\bf Sara BROFFERIO and Wolfgang WOESS}
\address{\parbox{1.4\linewidth}{Laboratoire de Math\'ematiques, 
Universit\'e Paris-Sud, B\^atiment 425,\\ F-91405 Orsay Cedex, France\\}}
\email{Sara.Brofferio@math.u-psud.fr}
\address{\parbox{1.4\linewidth}{Institut f\"ur Mathematik C, 
Technische Universit\"at Graz,\\
Steyrergasse 30, A-8010 Graz, Austria\\}}
\email{woess@TUGraz.at}
\date{January 1, 2005}
\thanks{Supported by European Commission, Marie Curie Fellowship
HPMF-CT-2002-02137}
\begin{abstract}
We determine all positive harmonic functions 
for a large class of ``semi-isotropic'' random walks on
the lamplighter group, i.e., the wreath
product $\Z_q \wr \Z$, where $q \ge 2$.
This is possible via the geometric realization of a Cayley graph 
of that group as the Diestel-Leader graph $\DL(q,q)$. More generally,
$\DL(q,r)$ ($q,r \ge 2$) is the horocyclic product of two homogeneous 
trees with respective degrees $q+1$ and $r+1$, and our result
applies to all $\DL$-graphs. This is based on a careful study
of the minimal harmonic functions for semi-isotropic walks on trees.
\end{abstract}
\maketitle

\markboth{\sf S. Brofferio and W. Woess}{\sf Positive harmonic functions}
\baselineskip 15pt

\section{Introduction}\label{intro}

Let $X$ be an infinite, connected, locally finite graph $X$ with root vertex $o$,
and $P$ the transition matrix $P = \bigl( p(x,y) \bigr)_{x,y \in X}$  of a random walk 
$(Z_n)_{n\ge 0}$ on $X$. That is, $Z_n \in X$ is the random position of the random walker at
time $n$, and $\Prb[Z_{n+1} = y \mid Z_n=x ] = p(x,y)\,$. 
The $n$-step transition probabilitiy
$p^{(n)}(x,y) = \Prb[Z_n=y \mid Z_0 = x]\,$, $x, y \in X\,$,
is the $(x,y)$-entry of the matrix power $P^n$, with $P^0=I$, the
identity matrix. The \emph{Green kernel} is
\begin{equation}\label{green}
G(x,y) = \sum_{n=0}^{\infty} p^{(n)}(x,y)\,.
\end{equation}
We suppose here that $P$ is \emph{irreducible:} 
$G(x,y) > 0$ for all $x,y \in X\,$, and \emph{transient:} $G(x,y) < \infty\,$.

A function $h: X \to \R$ is called \emph{harmonic,} or \emph{$P$-harmonic,}
if $h = Ph\,,$ where $Ph(x) = \sum_y p(x,y)\,h(y)\,$, and \emph{superharmonic}
if $h \ge Ph$.
A function $h \in \HH^+$, the cone of positive harmonic functions, 
is called \emph{minimal} if
$h(o) = 1\,$, and $h \ge h_1 \in \HH^+$ implies that $h_1/h$ is constant.
The minimal harmonic functions are the extreme points of the convex base
$\BB = \{ h \in \HH^+ : h(o)=1 \}$ of the cone $\HH^+$.  
Every positive harmonic function has a unique integral representation 
with respect to a Borel measure on the set of minimal ones, see 
{\sc Doob} \cite{Do}.

Positive harmonic functions for various classes of random walks (Markov chains)
have been a continuous subject of study since the 1950ies. One of the typical questions
is to determine and describe all positive harmonic functions in terms of a
geometric or algebraic structure of the underlying graph $X$, to which the
transition probabilities are assumed to be adapted. See the monograph by 
{\sc Woess} \cite{Wbook}, Ch. IV, for various results in this spirit.

In the present paper, we determine all minimal, and thereby also all positive
harmonic functions for a large class of random walks on the \emph{lamplighter
group} $\Z_q \wr \Z$. This is the \emph{wreath product} of the additive
group $\Z_q = \{ 0, \dots, q-1\}$ of integers modulo $q (\ge 2)$ with
the group $\Z$ of all integers. 

The lamplighter interpretation is as follows: $\Z$ represents an infinite
street, i.e. the graph with edges $[k,k+1]$, $k \in \Z$, with a lamp at the 
midpoint $k-\frac{1}{2}$ of each edge. Each lamp can have $q$ different states
in $\Z_q\,$; the state $0$ corresponds to the lamp being switched of. 
A lamplighter wanders along $\Z$ (from a point $k$ to $k\pm 1$), 
and at each step, he may change the state
of the lamp on the edge which he traverses. (Below, we will 
also allow bigger ``jumps'' along $\Z$ and changing the lamps
on more than one of the nearby edges).
For a corresponding random walk, the information that we have to keep track of
at each instant is the pair $(\eta,k)$, where $k$ is the current position
of the lamplighter, and $\eta$ is the current configuration of the 
states of the lamps.

We remark that for the usual construction of the wreath product, one thinks
of the lamps sitting at the points of $\Z$. For our purpose, it is more
convenient to have them (equivalently) sitting at the edges' midpoints,
i.e., the elements of $\Z - \frac{1}{2}$. In these terms, the formal 
construction of the lamplighter group is as follows.
Consider the group of all finitely supported \emph{configurations}
$
\CC = \{ \eta: \Z-\frac{1}{2} \to \Z_q \,,\;\supp(\eta) \;\ \mbox{finite}\}
$
with pointwise addition modulo $q$. Then $\Z$ acts on $\CC$
by translations $k \mapsto T_k:\CC \to \CC$ with 
$T_k\eta(m - \frac{1}{2}) = \eta(m-k-\frac{1}{2})$.
The resulting semidirect product $\Z \rightthreetimes \CC$ is 
$$
\Z_q \wr \Z = \{ (\eta,k) : \eta \in \CC\,,\;k \in \Z\} \,,
\quad \mbox{group operation}\quad
(\eta,k)(\eta',k') = (\eta + T_k\eta', k+k').
$$
The group identity is $o = (\0,0)$, where $\0$ is the zero configuration.
For two pairs $x=(\eta,k)$, $y=(\eta',k')$, we define the left and right 
\emph{flags} 
\begin{equation}\label{flags}
\begin{aligned}
\fl_1 &= \fl_1(x,y) = \min \{ k, k', m : \eta(m-\tfrac{1}{2}) \ne \eta(m-\tfrac{1}{2}) \}
\AND\\
\fl_2 &= \fl_2(x,y) = \max \{ k, k', m : \eta(m+\tfrac{1}{2}) \ne \eta(m+\tfrac{1}{2})\}.
\end{aligned}
\end{equation}
These are the left- and rightmost positions on $\Z$ which the lamplighter
is forced to visit if he starts at $k$ with configuration $\eta$ and
wants to reach $k'$ with configuration $\eta'$, when at each single step he
traverses a single edge and is allowed to change the state of the lamp on
that edge. We also define the corresponding increments
\begin{equation}\label{ups}
\up_1 = 
k - \fl_1 \AND 
\up_2 = 
\fl_2\ - k\,.
\end{equation}
We say that a random walk on $\Z_q \wr \Z$ is \emph{semi-isotropic} if 
the transition probability from $x=(\eta,k)$ to $y=(\eta',k')$ depends 
only on $\up_1$, $\up_2$ and $k-k'$.
Every random walk of this type is adapted to the group structure of
$\Z_q \wr \Z$. Indeed, $p(x,y) = \mu(x^{-1}y)$, where $\mu(x)=p(o,x)$,
a probability measure on $\Z_q \wr \Z$. 

A typical class of examples can be obtained 
as follows. For $m \in \Z \setminus \{0\}$, let $\mu_m$ be the probability 
measure associated with the random walk, where from position $k$ and 
configuration $\eta$, the lamplighter jumps to $k+m$ and switches each
of the lamps on the $|m-1|$ edges between $k$ and $k+m$ to a uniformly chosen
random state
(independendtly of each other), while leaving the other lamps unchanged.
Write $\mu_0=\de_o$, the point mass at the identity.
If $\wt\mu$ is any probability measure on $\Z$ then the probability
measure 
\begin{equation}\label{semi}
\mu = \sum_m \wt\mu(m)\,\mu_m
\end{equation}
gives rise to a semi-isotropic 
random walk. The latter is irreducible if and only if the random walk on $\Z$
induced by $\wt\mu$ is irreducible. Also, it is transient, since 
$\Z_q \wr \Z$ has exponential growth, see {\sc Varopoulos} \cite{Var}
or the exposition in \cite{Wbook}, Ch. I.

There are natural projections 
$\pi_1, \pi_2 : \Z_q \wr \Z \to \T_q$, where $\T_q$ is the homogenous
tree with degree $q+1$. Under each of the two projections, every semi-isotropic 
transition matrix $P$ projects to transition matrices $P_1$ and $P_2$
(respectively) on $\T_q$, which are also semi-isotropic in an
adequate sense.

Our results arise as a special case of a more general class of lamplighter 
type random walks, which -- as well as the ones discussed so far -- arise as
random walks on the \emph{Diestel-Leader graphs} $\DL(q,r)$, where
$q,r \ge 2$. In that description, the projections $\pi_1, \pi_2$ map
$\DL(q,r)$ onto $\T_q$ and $\T_r$, respectively. 
The details will be explained in \S 2. 

Our first result, Theorem \ref{main1}, states that every minimal
$P$-harmonic function $h$ is of the form $h(x)=h_1(\pi_1x)$ or
$h(x)=h_2(\pi_2x)$, where $h_i$ is a $P_i$-harmonic function on $\T_q$
($i=1,2$). This leads us to a careful study of all minimal harmonic functions 
for semi-isotropic random walks on $\T_q$. This is done under suitable moment 
conditions in Theorem \ref{main2}, on the basis of recent work of
{\sc Brofferio} \cite{Bro}. In Theorem \ref{main3}, we then describe all
minimal, and thereby all positive $P$-harmonic functions on $\DL(q,r)$.
The results and their proofs are a considerable extension as well as
simplification of those of {\sc Woess} \cite{Wo}, who 
only dealt with the nearest neighbour case. We also remark here that
another extension of \cite{Wo} is given by {\sc Brofferio and Woess}
\cite{BrWo}, who study only nearest neighbour random walks, but give
precise asymptotic estimates (in space) of the Green kernel, which
leads to a description of the full \emph{Martin compactification.}
The latter contains more analytical information than the one provided by
knowledge of the minimal
harmonic functions. However, it seems hard to extend the methods of \cite{BrWo} 
to general semi-isotropic random walks as considered in the present paper.

In concluding the introduction, let us remark that the first to show that
lamplighter groups are fascinating objects in the study of random
walks were {\sc Kaimanovich and Vershik} \cite{KaiVer}. 
By now, there is a considerable amount of literature on this
topic, regarding various issues. See e.g. 
{\sc Kaimanovich} \cite{Kai}, {\sc Lyons, Pemantle and Peres} 
\cite{LyoPemPer}, {\sc Erschler} \cite{Ers1}, \cite{Ers2}, 
{\sc Revelle} \cite{Rev1}, \cite{Rev2},
{\sc Bertacchi} \cite{Ber}, {\sc Grigorchuk and \.Zuk}
\cite{GriZuk}, {\sc Dicks and Schick} \cite{DiSc}, {\sc Bartholdi and Woess} 
\cite{BaWo}, {\sc Saloff-Coste and Pittet} \cite{PitSal1}, \cite{PitSal2}.

We also remark here that in part, our methods have their roots in the study of
random walks on the affine group over the real, resp. $p$-adic numbers,
see {\sc Elie} \cite{Eli}, {\sc Babillot, Bougerol and Elie} 
\cite{BaBoEl},  resp. {\sc Cartwright, Kaimanovich and Woess} \cite{CaKaWo}
and {\sc Brofferio} \cite{Bro}.

\section{Lamplighters and Diestel-Leader graphs}\label{geometry}

We now explain very briefly the structure of the DL-graphs and their 
relation with the wreath products $\Z_q \wr \Z$. See also \cite{Wo}, 
\cite{BaWo} and \cite{BrWo}; here we choose a different order of explanations,
which together with those of the latter papers may create a more complete
picture.

Consider the two-way-infinite path $\Z$ with lamps sitting at the midpoints 
of the edges, as above. However, we now think of a more general model,
where each lamp may be switched on in two different colours, green and red.
There are $q$ possible green states (intensities), encoded by 
$\Z_q$, and $r$ possible red states, encoded by $\Z_r$. In both cases,
the respective $0$ state means ``switched off''. Only finitely many lamps may be switched
on, and the rule is that all lamps on the left (towards $-\infty$) of the
lamplighter have to be in a green state, while all lamps on the right 
(towards $+\infty$) have to be in a red state. For each $k \in \Z$, let
$$
\CC_k = \{ \eta: \Z-\tfrac12 \to \Z_q \dot\cup \Z_r \mid \supp(\eta) 
\;\mbox{finite}\,,
\; \eta(m-\tfrac12) \in \Z_q \forall \ m \le k\,,
\; \eta(m+\tfrac12) \in \Z_r \forall \ m \ge k \}\,.
$$
The state space of our lamplighter walks, 
i.e., the vertex set of the $\DL$-graph, is the set
$X = \bigcup_{k \in \Z} \CC_k \times \{k\}\,$. Two pairs $(\eta,k)$ and
$(\eta',k')$ are neighbours if $k' = k\pm1$ and $\eta$ coincides with
$\eta'$ everywhere except at the midpoint of the edge between $k$ and $k'$.
This describes the Diestel-Leader graph $\DL(q,r)$.

In general, we consider $\Z_q$ and $\Z_r$ as being disjoint,
but when $q=r$, we need not distinguish between the two colours. In the latter
case, we can omit the index $k$ of $\CC_k$. It is obvious from this 
description that $\DL(q,q)$ is a Cayley graph of $\Z_q \wr \Z$, since 
the latter group acts on $\DL(q,q)$ transitively and fixed-point-freely 
by graph isometries.

We next want to achieve a geometric understanding of the graph structure.
If $\eta \in \CC_k$ then set
$$
\eta_k^-(m) = \eta(k-m-\tfrac12) \AND \eta_k^+(m) = \eta(k+m+\tfrac12)\,,
\quad m \ge 0\,.
$$
That is, we split $\eta$ at $k$ and consider the left ($\eta_k^-$) 
and right ($\eta_k^+$) halves
as elements of $\Sigma_q$ and $\Sigma_r$, respectively, where
$$
\Sigma_q = \bigl\{ \bigl(\sigma(m)\bigr)_{m \ge 0} : \sigma(\cdot) \in \Z_q\,,\;
\supp(\sigma) \;\mbox{finite} \bigr\}\,.
$$
(In \cite{Wo}, these sequences were indexed over the non-positive integers.)
The set $\Sigma_q \times \Z$ is in one-to-one correspondence with
the homogeneous tree $\T_q$ of degree $q+1$. Seen as a vertex of that tree, each
element $(\sigma, k)$ has a unique \emph{predecessor} $(\sigma',k-1)$,
where $\sigma'(m) = \sigma(m+1)$ for $m \ge 0$, that is, $\sigma'$ is obtained
by deleting $\sigma(0)$ from $\sigma$. This describes neighbourhood, so
that $(\sigma,k)$ has precisely $q$ \emph{successors} - the elements that
have $(\sigma,k)$ as their predecessor. 

We can draw this tree in horocyclic layers. The $k$-th \emph{horocycle} is
$H_k = \Sigma_q \times \{ k \}$, and for every element of $H_k$, its predecessor
lies on $H_{k-1}$, while its successors lie on $H_{k+1}$.
See Figure 1a below, and Figure 1 in \cite{Wo} for a more detailed picture. 
The closure (as a partial order) of the predecessor relation is the 
\emph{ancestor} relation. 
Every pair of vertices $x_1=(\sigma, k)$ and $y_1=(\sigma', k')$ has an infimum,
i.e., a maximal common ancestor. Formally, this is 
$x_1 \cf y_1 = (\bar \sigma, \bar k)$, where
$$
\bar k = \max \bigl\{ \ell \le \min\{k,k'\} : 
\sigma|_{[k-\ell,\,,\infty)} \equiv \sigma'|_{[k'-\ell,\,,\infty)} \bigr\}
\AND \bar\sigma(m) = \sigma(k - \bar k + m)\,,\; m\ge 0\,.
$$
In this notation, setting $\up(x_1,y_1) = k-\bar k$ (whence
$\up(y_1,x_1) = k' - \bar k$), the graph distance in $\T_q$
is $d(x_1,y_1) = \up(x_1,y_1) + \up(y_1,x_1)$, and $\up(y_1,x_1) = 
d(x_1,x_1 \cf y_1)$. Also, $\up(y_1,x_1) - \up(x_1,y_1) = \hor(y_1)-\hor(x_1)$,
where $\hor(x_1) = k$ if $x \in H_k$. See Figure 1b.

$$
\beginpicture
\setcoordinatesystem units <3mm,3.8mm> 

\setplotarea x from -4 to 30, y from -3.8 to 6.4
\arrow <5pt> [.2,.67] from 4 4 to 1 7
\put{$\omega_1$} [rb] at 1.2 7.2


\plot -4 -4       4 4         /         
\plot 4 4         12 -4          /      
\plot -2 -2       -2.95 -4 /            
\plot -.5 -2      -1.9 -4     /         
\plot -.5 -2      -.85 -4   /           
\plot 1 -2        .2 -4      /          
\plot 1 -2        1.25  -4     /        
\plot 2.5 -2      2.3  -4     /         
\plot 2.5 -2      3.35 -4      /        
\plot 5.5 -2      4.65  -4    /         
\plot 5.5  -2     5.7  -4      /        
\plot 7  -2       6.75   -4   /         
\plot 7 -2        7.8  -4      /        
\plot 8.5  -2     8.85  -4    /         
\plot 8.5  -2     9.9  -4      /        
\plot 10  -2      10.95  -4   /         
\plot 0  0        -.5  -2      /        
\plot 2 0         1 -2     /            
\plot 2 0         2.5   -2     /        
\plot 6 0         5.5    -2    /        
\plot 6 0         7 -2         /        
\plot 8 0         8.5 -2       /        
\plot 2 2         2 0         /         
\plot 6 2         6 0         /         



\put {$\vdots$} at 4 -5.2

\put {\it Figure 1a} at 4 -7.5


\put{$x_1$} [lb] at  32.15 -2.8
\put{$y_1$} [rb] at 25.85 -0.8
\put{$x_1 \cf y_1$} [lb] at 28.15 1.2

\arrow <5pt> [.2,.67] from 28 1 to 24 5 
\plot 32 -3  28 1  26 -1 /
\multiput {$\scs\bullet$} at 32 -3  28 1  26 -1 /

\setdots<4pt>
\plot 22 1  36 1 /
\plot 22 -1 26 -1 /
\plot 32 -3  36 -3 /

\put{$\up(y_1,x_1)\!\left\{\rule[-2.8mm]{0mm}{0mm}\right.$}[rb] at 22 -1 
\put{$\!\!\left.\rule[9mm]{0mm}{0mm}\right\}\!\up(x_1,y_1)$}[lb] at 36 -3

\put {\it Figure 1b} at 28 -7.5
\put{$\;$} at 28 -8.5

\endpicture
$$

Given $\DL = \DL(q,r)$, we write $\T^1 = \T_q$ and $\T^2 = \T_r$ and define the
projections $\pi_i: \DL \to \T_i$ as follows. If $x = (\eta,k)$
then $\pi_1x = x_1 := (\eta_k^-,k)$ and $\pi_2x = x_2 :=(\eta_k^+,-k)$.
Each of these mappings is a neighbourhood-preserving surjection of
$\DL$ onto the respective tree. Conversely, starting with the trees, 
$$
\DL(q,r) = \{ x=x_1x_2 : x_i \in \T^i\,,\; \hor(x_1) + \hor(x_2) = 0 \}
$$
with neighbourhood given by $x_1x_2 \sim y_1y_2 \iff x_1 \sim y_1$ and
$x_2 \sim y_2$. Again, a detailed geometric picture can be found in \cite{Wo}, 
Figure 2. In that reference (as well as in \cite{BaWo} and \cite{BrWo}),
the explanation follows the reversed order, starting with the geometric 
description. We recall that when $x=x_1x_2, y=y_1y_2 \in \DL$ then
\begin{equation}\label{upup}
\up(x_1,y_1) + \up(x_2,y_2) = \up(y_1,x_1) + \up(y_2,x_2)\,,
\end{equation}
and their distance is,  by \cite{Ber},
\begin{equation}\label{dist}
d(x,y) = d(x_1,y_1)+d(x_2,y_2) - |\hor(x_1)-\hor(y_1)|\,.
\end{equation}
With the 2-colour-lamplighter interpretation, the left and right flags
can be defined as in \eqref{flags}, and if $x=(\eta,k)$, $y=(\eta',k')$
then
\begin{equation}\label{flags2}
\fl_i= \hor(x_i \cf y_i)\,,\; i=1,2\,.
\end{equation}
Every $\DL$-graph is vertex-transitive (its isometry group acts transitively
on the vertex set), but only when $r=q$ it is a Cayley graph of a
finitely generated group. 

\section{Positive harmonic functions on $\DL$-graphs}

A random walk, resp. its transition matrix $P_i$ on $\T^i$ ($i=1,2$)
is called \emph{semi-isotropic}, if $p_i(x_i,y_i)$ depends only on 
$\up(x_i,y_i)$ and $\up(y_i,x_i)$, or equivalently, only on
$d(y_i,x_i)$ and $\hor(y_i)-\hor(x_i)$. (Recall that a random walk is
called isotropic if $p_i(x_i,y_i)$ depends only on the distance.)

Anologously, we call a random walk, resp. its transition matrix $P$ on $\DL$
\emph{semi-isotropic}, if $p(x,y)$ depends only on the four numbers
$\up(x_1,y_1)$, $\up(y_1,x_1)$, $\up(x_2,y_2)$ and $\up(y_2,x_2)$,
which must satisfy \eqref{upup}. That is, there is a probability measure $\mm$ on the set 
$\bigl\{(k_1,l_1,k_2,l_2) \in 
\N_0^4 : k_1+k_2=l_1+l_2 \bigr\}$ (where $\N_0$ denotes the non-negative
integers) such that
\begin{equation}\label{isotropic}
p(x,y) = \mm\bigl(\up(x_1,y_1),\up(y_1,x_1),\up(x_2,y_2),\up(y_2,x_2)\bigr)
\quad\
\end{equation}
for all $x=x_1x_2\,,\; y=y_1y_2 \in \DL(q,r)\,$. In terms of the lamplighter, 
the transition probability from $(\eta,k)$ to $(\eta',k')$ depends
only on the distances of $k$ and $k'$ to the left and right flags \eqref{flags},
\eqref{flags2}.

The projection $P_1$ of $P$ on $\T^1$ is given by
$$
p_1(x_1,y_1) = \sum_{y_2 : y_1y_2 \in \DL} p(x_1x_2,y_1y_2)\,,
$$
which is independent of the specific choice of $x_2 \in \T^2$ such that 
$x_1x_2 \in \DL$. The projection $P_2$ on $\T^2$ is analogous. Both are
semi-isotropic along with $P$.

\begin{lem}\label{transient} Every semi-isotropic, irreducible random
$P$ walk on $\DL(q,r)$ is transient, and its projections $P_i$ on 
$\T^i$ ($i=1,2$) are also transient.
\end{lem}

\begin{proof} This follows from Theorem 5.13 in \cite{Wbook}, 
since $\DL(q,r)$ as well as the $\T^i$ are vertex-transitive graphs with
exponential growth.
\end{proof}

Let $f$ be a function defined on $\DL(q,r)$. When we say that 
\emph{$f$ depends only on $x_1$}, resp. \emph{$f$ depends only on $x_2$},
then this means
$$
f(x_1x_2) = f(x_1y_2) \quad\forall\ x_1x_2, x_1y_2 \in \DL\,,
\quad \mbox{resp.} \quad 
f(x_1x_2) = f(y_1x_2) \quad\forall\ x_1x_2, y_1x_2 \in \DL\,.
$$
In the first case, we can write $f(x_1x_2) = f_1(x_1)$, where
$f_1$ is a function on $\T^1$, and in the second case,
$f(x_1x_2) = f_2(x_2)$, where $f_2$ is a function on $\T^2$.
Note that when both conditions hold, then this does not mean that
$f$ is constant, but that $f(x_1x_2) = \wt f\bigl(\hor(x_1)\bigr)$,
where $\wt f$ is a function on $\Z$. The following is an obvious
exercise.

\begin{lem}\label{project-harmonic} 
A function $h_1$ is $P_1$-harmonic on $\T^1$ if and only if 
$h(x_1x_2) = h_1(x_1)$ is $P$-harmonic on $\DL$.
\end{lem}

Here is the first main main result, along with a surprisingly simple proof.

\begin{thm}\label{main1}
Suppose that $P$ is semi-isotropic and irreducible on $\DL$.
Then the following statements hold.

\smallskip

{\rm(a)} Every minimal $P$-harmonic function on $\DL(q,r)$
is of the form $h(x_1x_2) = h_1(x_1)$ or $h(x_1x_2) = h_2(x_2)$,
where $h_i$ is a minimal $P_i$-harmonic function on $\T^i$, $i=1,2$
(respectively).

\smallskip

{\rm (b)} If $h$ is a positive $P$-harmonic function on $\DL(q,r)$,
then there are non-negative $P_i$-harmonic functions $h_i$ on
$\T^i$ ($i=1,2$) such that
$$
h(x_1x_2) = h_1(x_1) + h_2(x_2) \quad \forall\ x_1x_2 \in \DL(q,r)\,.
$$
\end{thm}

We shall appeal to \emph{Martin boundary} theory for
Markov chains, see {\sc Doob} \cite{Do}, {\sc Hunt} \cite{Hu}, 
or the excellent introduction by {\sc Dynkin} \cite{Dy}. 
If $G(\cdot,\cdot)$ is the Green kernel of any transient, irreducible
Markov chain with countable state space $X$, then the
\emph{Martin kernel} is
$$
K(x,y) = G(x,y)/G(o,y)\,,\quad x, y \in X\,,
$$
where $o \in X$ is a reference point. The \emph{Martin compactification}
is the smallest compact Hausdorff space $\wh X$ containing $X$ as a dense, 
discrete subset, such that for each $x \in X$, the function $K(x,\cdot)$ 
extends continuously to $\wh X$. The extendend kernel on $X \times \wh X$
is also denoted $K(\cdot,\cdot)$. The \emph{Martin boundary} is 
$\MM = \wh X \setminus X$. A basic result of the theory is that every 
minimal harmonic function is of the form $K(\cdot, \zeta)$, where 
$\zeta \in \MM$, and that the \emph{minimal Martin boundary}
$\MM_{\min}$, consisting of all $\zeta \in \MM$ for which $K(\cdot,\zeta)$ 
is minimal harmonic, is a Borel subset of $\MM$. The \emph{Poisson-Martin
representation theorem} says that for every positive harmonic function $h$
there is a unique Borel measure $\nu^h$ on $\MM$ such that
\begin{equation}\label{Poisson-Martin}
\nu^h(\MM \setminus \MM_{\min}) = 0 \AND 
h(x) = \int_{\MM} K(x,\cdot)\, d\nu^h \quad \forall\ x \in X\,.
\end{equation}

\begin{proof}[Proof of Theorem \ref{main1}]  
Recall that our ``root'' (origin) of $\DL$ is $o = (\0,0)$. Let $h = K(\cdot,\zeta)$
be a minimal harmonic function. Then there is a sequence 
$y^{(n)} = y^{(n)}_1y^{(n)}_2\in \DL$ such that 
$$
K(x,y^{(n)}) \to h(x) \quad \forall\ x \in \DL\,.
$$
Then, by \eqref{upup} and \eqref{dist}, at least one of the sequences
$\bigl(\up(o_i,y_i^{(n)})\bigr)_n$ is unbounded. (Recall that for any
$x \in \DL$, we write $x_i=\pi_ix$ for its projection on the tree $\T^i$.)
Thus, passing to a subsequence, we assume that
$\up(o_1,y_1^{(n)}) \to \infty$ (first case). 
We claim that in this case, $h(x)$ depends only on $x_2$. 
Indeed, fix $x = x_1x_2 \in \T^2$ and let $v_1$ lie on the same horocycle of 
$\T^1$ as $x_1$, so that also $v = v_1x_2 \in \DL$. 
Let $k = \up(x_1,v_1) = \up(v_1,x_1)$. By assumption,
$\up(x_1,y_1^{(n)}) \to \infty$ as $n \to \infty$.
Therefore there is $n(k)$ such that $\up(x_1,y_1^{(n)}) > k$ for all 
$n \ge n(k)$. This implies
\begin{equation}\label{up-up}
\up(x_1,y_1^{(n)}) = \up(v_1,y_1^{(n)}) \AND
\up(y_1^{(n)},x_1) = \up(y_1^{(n)},v_1) \quad \forall\ n \ge n(k)\,,
\end{equation}
see Figure 2. 
$$
\beginpicture
\setcoordinatesystem units <3.2mm,3.7mm> 

\setplotarea x from 0 to 8, y from 0 to 8

\put{$x_1$} [lb] at  8.15 0.2
\put{$v_1$} [rb] at -.15  0.2
\put{$y_1^{(n)}$} [rb] at 0.7 4.7

\arrow <5pt> [.2,.67] from 4 4 to 0 8 
\plot 8 0  4 4  0 0 /
\plot 2 6  0.5 4.5  /

\multiput {$\scs \bullet$} at 8 0  0 0  2 6  4 4  0.5 4.5   /

\setdots<4pt>
\plot -2 0  10 0  /



\endpicture
$$

\smallskip

\begin{center}
\emph{Figure 2}
\end{center}

\vspace{.4cm}

Formula \eqref{isotropic} implies that also the Green kernel is semi-isotropic,
and consequently
$$
K(x,y^{(n)}) 
= \frac{G(x_1x_2,y_1^{(n)}y_2^{(n)})}{G(o_1o_2,y_1^{(n)}y_2^{(n)})} = 
\frac{G(v_1x_2,y_1^{(n)}y_2^{(n)})}{G(o_1o_2,y_1^{(n)}y_2^{(n)})}
= K(v,y^{(n)})\,.
$$
for all $n \ge n(k)$. Letting $n \to \infty$, we obtain
$h(x) = h(v)$. Lemma \ref{project-harmonic} yields that 
$h(x) = h_2(x_2)$, where $h_2$ is a $P_2$-harmonic function on 
$\T^2$. Minimality of $h$ as a $P$-harmonic function implies 
minimality of $h_2$ as a $P_2$-harmonic function. (The converse is in 
general not true.)

By exchanging the roles of the two trees, we see that if
$\up(o_2,y_2^{(n)}) \to \infty$ (second case) then $h(\cdot)$ 
depends only on $x_2$, and $h(x) = h_1(x_1)$, where $h_1$ is
a minimal $P_1$-harmonic function on $\T^1$.

This proves (a). 
To see (b), let 
$$
\MM_i(P) = \{ \zeta \in \MM(P) : K(\cdot, \zeta) \;\mbox{is minimal
and depends only on}\; x_i \}\,,\quad i=1,2\,.
$$
Then $\MM_{\min}(P) = \MM_1(P) \cup \MM_2(P)$. (The two pieces are
not necessarily disjoint, see below.) The topology of the Martin boundary
is the one of pointwise convergence of the Martin kernels, and $\MM_{\min}$
is a Borel subset. Since the set of $\SH_i(P)$ of all positive $P$-superharmonic
functions that depend only on $x_i$ is closed with respect to pointwise
convergence\footnote{When $P$ does not have finite range, we cannot
take positive $P$-harmonic functions here, since a pointwise limit of such 
functions is not necessarily harmonic, while being superharmonic by Fatou's
theorem}, $\MM_i(P) = \MM_{\min}(P) \cap \SH_i(P)$ is a Borel
subset of the Martin boundary. Thus, if $\nu^h$ is as in \eqref{Poisson-Martin},
then $h = h_1 + h_2$, where
\begin{equation}\label{h-sum}
h_1 = \int_{\MM_1(P)} K(\cdot,\zeta)\,d\nu^h(\zeta) \AND
h_2 = \int_{\MM_2(P) \setminus \MM_1(P)} K(\cdot,\zeta)\,d\nu^h(\zeta)\,.
\end{equation}
By (a), $h_i$ is $P_i$-harmonic for $i=1,2$.
\end{proof}

\section{Semi-isotropic random walks on a homogeneous tree}\label{tree}

In view of Theorem \ref{main1}, our next aim is to determine those
minimal harmonic functions for $P_i$ on $\T^i$ ($i=1,2$) which lift to
minimal $P$-harmonic functions on $\DL(q,r)$. Note that when $h_i$ is minimal 
harmonic for $P_i$ on $\T^i$  then $h(x_1x_2) = h_i(x_i)$ is not necessarily
minimal for $P$ on $\DL$ (while the converse \emph{is} true), compare with
{\sc Woess} \cite{Wo}. 

In the following we shall omit the subscripts $_i$ for elements of $\T^i$. 
Thus, \emph{in the present section,
$P$ denotes an irreducible, semi-isotropic random walk on $\T= \T_q$.}

We recall the construction of the geometric boundary $\bd \T$. A 
\emph{geodesic ray} is a one-sided infinite path in $\T$ without 
repeated vertices. A boundary point (end) is an equivalence class of rays,
where two rays are equivalent if they differ only by finitely many
initial points. If $x \in \T$ and $w \in \wh \T = \T \cup \bd \T$ 
then there is a unique geodesic path from $x$ to $w$ (if $w \in \T$), 
resp. ray representing $w$ (if $w \in \bd \T$) starting at $x$, denoted
by $\geo{x\,w}$. Analogously, if $\xi, \eta$ are distinct ends,
then there is a unique two-sided infinite geodesic path $\geo{\xi\, \eta}$
whose two ``halves'' (when split at any of its vertices) represent $\xi$
and $\eta$.

The \emph{confluent} $c(v,w)$ of $v, w \in \wh \T$ is the last common
vertex on the geodesics from the origin $o \in \T$ to $v$ and $w$, respectively.
Writing $|x| = d(x,o)$ for $x \in \T$, we can equip $\wh \T$ 
with the ultrametric $\theta(v,w) = q^{-|c(v,w)|}$, if $v \ne w$
(and $\theta(v,v)=0$). Thus, $\wh T^1$ is a compact space with $\T$ discrete,
open and dense. 

Let $\om$ be the end of $\T$ represented by the ray 
$[o=o_0, o_1, o_2, \dots]$ where each $o_k$ is the predecessor
of $o_{k-1}$, and the root $o$ is on the horocycle $H_0$. 
In the correspondence between vertices of $\T$ and sequences described in 
\S 2, we have $o_k = (\0,-k)$.
Then $\om=\om_1$ is the end of $\T=\T^1$ located at the top of the picture
in Figure 1a, while the ends in $\bd^*\T = \bd\T \setminus \{\om\}$ 
are located at the bottom of that picture. Each $\xi\in \bd\T$
corresponds to a two-sided infinite sequence in $\Z_q^{\Z}$
for which there is $\bar m \in \Z$ such that $\xi(n) = 0$ for all
$n \le \bar m$. Given $x \in \T$ and $\xi \in \bd^*\T$,
we can define $x \cf \xi$ and $\up(x,\xi)$ in the same way as
above (see Figure 1b).

Now consider the Martin kernel $K(x,y) = G(x,y)/K(o,y)$
associated with $P$ on $\T$. In the case when $P$ has 
\emph{bounded range},
i.e., $p(x,y) = 0$ when $d(x,y) > R$ ($R < \infty$) 
then it is easy to find the minimal $P$-harmonic functions as a
consequence of a general result of {\sc Picardello and Woess} \cite{PiWo}
regarding transient, bounded range random walks on arbitrary trees.

\begin{pro}\label{bounded-range}
If $P$ is irreducible and has bounded range then the associated 
Martin compactification is the end compactification $\wh\T$, and
each extended kernel $K(\cdot,\xi)$, $\xi \in \bd \T$ is
a minimal harmonic function for $P$.
\end{pro}


We want to extend the resulting description of the minimal Martin boundary 
in two ways. First, we go beyond bounded range, and second,
we shall describe the resulting minimal Martin kernels in more computational
detail.

\medskip

{\bf Invariant measures on the boundary.} 
The transition probabilities of $P$ are invariant under the locally compact, 
totally disconnected group $\Ga=\AFF(\T)$ of all isometries of $\T$ that fix $\om$.
This group acts transitively on $\T$, so that we can interpret our random 
walk as a random walk on that group via the construction described, e.g.,
in {\sc Woess} \cite{Wo1}, \S 3. Namely, normalize the left Haar measure
$dg$ on $\Ga$ such that the stabilizer $\Ga_{\! o}$ of $o$ in $\Ga$ (which
is an open-compact subgroup) has measure $1$. Define a probability measure 
$\mu$ on $\Ga$ by 
\begin{equation}\label{mu}
\mu(dg) = p(o,go)\,dg\,.
\end{equation}
Let $X_n$, $n \ge 1$, be a sequence of i.i.d. $\Ga$-valued random variables
with common distribution $\mu$. Consider the \emph{right random walk} 
$$
R_n = X \cdots X_n
$$ 
on $\Ga$. Then, given $g \in \Ga$,
the sequence $g\,R_no$ is (a model of) the random walk on $\T$ with transition
matrix $P$ and starting point $x=go$. In particular, all results 
of {\sc Cartwright, Kaimanovich and Woess} \cite{CaKaWo} and 
{\sc Brofferio} \cite{Bro} regarding random walks on $\Ga$ apply here.
Since the action of $\Ga$ extends to $\bd^*\T$, we can convolve
$\mu$ with any (Radon) measure $\nu$ on $\bd^*\T$. If $E \subset \bd^*\T$
is a Borel set, then
$$
\mu*\nu(E) = \int_{\Ga} \nu(g^{-1}E)\,\mu(dg)\,.
$$
We are looking for an \emph{invariant} measure $\nu$, satisfying $\mu*\nu=\nu$.
It will serve to describe the minimal harmonic functions.
For $x\in \T\,$, $\xi \in \bd^*\T$ and $k, r \ge 0$, let 
\begin{equation}\label{Om_k}
\begin{gathered} 
T_{k,r}(x) = \{ y \in \T : \up(x,y)=k\,,\;\up(y,x)=r\} \,,\\
\Om_k(x) = \{ \eta \in \bd^*\T : \up(x,\eta)=k \}\AND
\Om(\xi) = \Om_0(x_0)\,,
\end{gathered}
\end{equation}
where $x_0$ is the element with $\hor(x_0)=0$ on the geodesic 
$\geo{\om\,\xi}$. If $x=o$ then we just write $\Om_k\,$ and $T_{k,r}$.
The sets $\Om_0(x)$,  $x \in \T$, are open-compact and generate the
topology of $\bd^*\T$. 

\begin{lem}\label{equidist}
If $\mu$ is defined by \eqref{mu} and $\nu$ is $\mu$-invariant on 
$\bd^*\T$, then $\nu$ is equidistributed on each set $\Om_k$, that is, 
for each $l \ge 1$ and $y \in T_{k,r}$ with $r \ge 1$, 
$$
\nu\bigl(\Om_0(y)\bigr) = \nu(\Om_k)/|T_{k,r}|\,.
$$
\end{lem}

\begin{proof} By \eqref{mu}, the measure $\mu$ is invariant under the 
stabilizer of $o\,$: if $g \in \Ga_{\! o}$ then 
$\de_g*\mu = \mu$. Therefore also $\de_g*\nu = \nu$. If 
$y, z \in T_{k,r}$
then there is $g \in \Ga_{\! o}$ such that $gy=z$. Thus
$\nu\bigl(\Om_0(z)\bigr) = \de_g*\nu\bigl(\Om_0(z)\bigr) = 
\nu\bigl(g^{-1}\Om_0(z)\bigr) = \nu\bigl(\Om_0(y)\bigr)$.
\end{proof}

The following lemma (due to Donald Cartwright)
is now the result of a straightforward computation
of the numbers \mbox{$\mu*\nu(\Omega_j)$}, $j \ge 0$.

\begin{lem}\label{a_k}
If $\mu$  and $\nu$ are as in Lemma \ref{equidist}, and
$$
\mu_{k,r} = \sum_{x \in T_{k,r}} p(o,x)\,,
$$
then the measure $\nu$ is determined
by the values $a_j = \nu(\Om_j)$, $j \ge 0$. 
The latter must satisfy the following equations.
$$
a_0
= \sum_{k=1}^{\infty} \sum_{r=1}^{\infty} \frac{a_r}{(q-1)q^{k-1}} \mu_{k,r}
+ \sum_{k=1}^{\infty} \frac{a_0}{q^k} \mu_{k,0}
+ \sum_{r=0}^{\infty} \Bigl(a_0+\dots+a_r\Bigr) \mu_{0,r}
$$
and for each $j \ge 1$,
$$
\begin{aligned}
a_j
&= \sum_{k=0}^{j-1} \sum_{r=0}^{\infty} a_{j+r-k} \mu_{k,r}
+ \sum_{k=j+1}^{\infty} \sum_{r=1}^{\infty} \frac{a_r}{q^{k-j}}\mu_{k,r}\\
&\quad + \frac{q-1}{q} \sum_{k=j}^{\infty} \frac{a_0}{q^{k-j}} \mu_{k,0}
+ \sum_{r=1}^{\infty} \Bigl(a_0+\dots+a_{r-1}+\frac{q-2}{q-1}a_r\Bigl) 
\mu_{j,r}\,.
\end{aligned}
$$
In particular, if $p(x,y)=0$ whenever $d(x,y) > N$ (bounded range), 
then
\begin{equation}\label{bounded}
a_j=\sum_{n=-N}^N a_{j+n}\wt\mu(n) \quad \text{for all}\quad j > N\,,
\quad\text{where}\quad\wt\mu(n) = \sum_{r-k=n} \mu_{k,r}\,\,.
\end{equation}
\end{lem}

Next, we look for sufficient conditions that guarantee the existence
(and uniqueness) of a solution of the above system for the $a_j$, or
equivalently, of a $\mu$-invariant measure $\nu$ on $\bd^*\T$.

Our basic requirement is that $P$ has \emph{finite first moment} $m(P)$, where
more generally the moment of order $t>0$ is defined as
\begin{equation}\label{moment}
m_t(P) = \sum_{x} d(o,x)^t\,p(o,x)
\end{equation}
(This is a generic definition, whenever we have a transition matrix 
and a metric.)

We can consider the projection $\wt P$ of $P$ onto $\Z$, where for 
arbitrary $x \in \T$ with $\hor(x) = k$, 
\begin{equation}\label{project-Z}
\wt p(k,l) = \sum_{y:\hor(y) = l} p(x,y) = \wt \mu(l-k)\,,
\end{equation}
with $\wt\mu(n)$ as in \eqref{bounded}. 
This defines 
an irreducible, translation-invariant random walk on $\Z$ which 
also has finite first moment, so that we can define its \emph{drift}
\begin{equation}\label{drift}
\al(P) = \al(\wt P) = \sum_{n \in \Z} n\,\wt p(0,n)\,.
\end{equation}

\begin{pro}\label{invariant} Suppose that $P$ is semi-isotropic and 
irreducible on $\T = \T_q$.\\[4pt]
{\rm (a)} If $m(P) < \infty$ and $\al(P) > 0$
then up to constant multiples, there is a unique $\mu$-invariant measure 
$\nu$ on $\bd^*\T$, and its total mass $\nu(\bd^*\T)$ is finite.\\[4pt]
{\rm (b)} If $m_{2+\ep}(P) < \infty$ (where $\ep > 0$) and $\al(P) = 0$
then up to constant multiples, there is a unique $\mu$-invariant measure 
$\nu$ on $\bd^*\T$, and $\nu(\bd^*\T)=\infty\,$.\\[4pt]
\indent
In both cases, $\nu$ is supported by the whole of $\bd^*\T$.
\end{pro}

\begin{proof}
(a) In this case, it is well known \cite{CaKaWo} that the random 
walk $Z_n=R_no$ starting at $o$ converges in the topology of $\wt \T$ to a 
$\bd^*\T$-valued random variable $Z_{\infty}$, and $\nu$ is its 
distribution. Irreducibility of $P$ implies that $\nu$ is supported by
the whole boundary. Uniqueness follows from the fact that 
$R_n \xi \to Z_{\infty}$ almost surely for every
$\xi \in \bd^*\T$, see \cite{CaKaWo}.

\smallskip

(b) This is proved in \cite{Bro}, Prop. 2.4.
\end{proof}

The case $\al(P) < 0$ is different. For general, semi-isotropic $P$,
consider the function 
\begin{equation}\label{phi}
\varphi(c) = \sum_{m \in \Z} \wt\mu(m)\,e^{cm}.
\end{equation}
It is well known that irreducibility implies that 
$\varphi: \R \to (0\,,\,\infty]$ is  convex (strictly convex where it is
finite) and that $\lim_{c \to \pm \infty} \varphi(c) = \infty\,$. 
Thus, there are at
most two solutions to the equation $\varphi(c)=1$. We have $\varphi(0)=1$,
and $\varphi'(0) = \al(P)$, if the derivative exists. In particular, 
if $\al(P) < 0$ and $\varphi(c_0) = 1$ for some $c_0 \ne 0$, then it must
be $c_0 > 0$; a sufficient condition for the existence of such a value
$c_0$ is that $\limsup_{n \to \infty} \wt\mu(n)^{1/n} = 0$. 
Given $c_0$, we define a new transition matrix $P^{\sharp}$ and associated
probability measure $\mu^{\sharp}$ on $\Ga$ by
\begin{equation}\label{conjugate}
p^{\sharp}(x,y) = p(x,y) e^{c_0(\hor(y)-\hor(x))} \AND
\mu^{\sharp}(dg) = e^{c_0\hor(go)}\mu(dg)
\end{equation}
$P$ is stochastic (since $\varphi(c_0)=1$), irreducible, and inherits 
semi-isotropy from $P$. 

\begin{pro}\label{invariant-neg}
If there is $c_0 \in \R$ such that
\begin{equation}\label{hor-moment}
c_0 > 0\,,\quad \varphi(c_0)=1\,,
\AND\sum_{x_i} d(o_i,x_i)\,p(o_i,x_i)\,e^{c_0\hor(x_i)} < \infty\,,
\end{equation}
then for $\mu^{\sharp}$ as in \eqref{conjugate}, 
up to constant multiples, there is a unique 
$\mu^{\sharp}$-invariant measure 
$\nu$ on $\bd^*\T$. Its support is $\bd^*\T$,
and $\nu(\bd^*\T)<\infty\,$.
\end{pro}
\begin{proof}
First of all, note that here we did not require existence of the
first moment. However, if $m(P) < \infty$, then \eqref{hor-moment}
implies $\al(P) < 0$ because of the shape of the function $\varphi(c)$.

Consider  $P^{\sharp}$ on $\T$, defined in \eqref{conjugate}.
Condition \eqref{hor-moment} says that
$m(P^{\sharp})<\infty$, whence $\al(P^{\sharp})$ is finite. If 
$\al(P^{\sharp})$ were non-negative, then we could not have 
$\varphi(c) = 1$ for any $c < c_0$, contradicting the fact that 
$\varphi(0)=1$.
Therefore, $\al(P^{\sharp}) > 0$, and we can apply Proposition 
\ref{invariant}(a) to
$P^{\sharp}$ and the associated probability measure $\mu^{\sharp}$
on $\Ga$. We find the unique $\mu^{\sharp}$-invariant probability
measure $\nu$ on $\bd^*\T$.
\end{proof}

\begin{rmk}\label{remark}
In order to compute the coefficients $a_j$ associated with $\nu$ via
Lemma \ref{a_k}, one has to replace the numbers $\mu_{k,r}$
with $\mu_{k,r}^{\sharp} = \mu_{k,r} \,e^{c_0(r-k)}$. 
\end{rmk} 

\medskip

{\bf Harmonic measures and Radon-Nikodym derivatives.} 

\smallskip

\underline{\sl A. Non-negative drift.} Suppose that $m(P) < \infty$ and 
$\al(P) > 0$, or that $m_{2+\ep}(P) < \infty$ and 
$\al(P) = 0$. Let $\nu$ be the invariant measure 
(probability measure in the first case) according to Proposition
\ref{invariant}. Define
$$
\nu_{x} = \de_g*\nu\,,\quad\text{where}\quad x \in \T
\AND g \in \Ga\;\;\text{with}\;\; go=x\,,
$$
that is, $\nu_{x}(E) = \nu(g^{-1}E)$ for Borel sets $E \subset \bd^*\T$.
Since $\nu$ is $\Ga_{\! o}$-invariant, the measure $\nu_{x}$
is independent of the specific choice of $g$ with $go=x$.
We observe that when $\al(P) > 0$ then $\nu_{x}(E)$ is the probability 
that the random walk governed by $P$ and starting at $x$ converges 
to a point in $E$. We have
\begin{equation}\label{harm-meas}
\nu_x = \sum_{y} p(x,y) \, \nu_{y} \quad \text{for every}\; x \in \T\,.
\end{equation}
Therefore, we call the family of measures $(\nu_{x})_{x \in \T}$
\emph{harmonic measures.}

\smallskip

\underline{\sl B. Negative drift.} Supose that \eqref{hor-moment} holds, and
let $\nu$ be the unique $\mu^{\sharp}$-invariant probability
measure on $\bd^*\T$ according to Proposition \ref{invariant-neg}.
This time, define
$$
\nu_{x} = e^{c_0\hor(x)}\, \de_g*\nu\,,\quad\text{where}\quad x \in \T
\AND g \in \Ga\;\;\text{with}\;\; go=x\,.
$$ 
Then again, the measures $(\nu_{x})_{x \in \T}$ satisfy \eqref{harm-meas}.

\begin{pro}\label{derivatives}
Under the moment conditions of \eqref{invariant} and
\eqref{hor-moment}, respectively, the harmonic measures are mutually
absolutely continuous. For each $x \in \T$, the Radon-Nikodym derivative
$d\nu_{x}/d\nu_{o}$ has a continuous realization, which we denote 
by $K(x,\cdot)$. 
\\[4pt]
If  $k,l$ are such that $\Om_k(o) \cap \Om_l(x) \ne \emptyset$, then
$$
K(x,\xi) = \frac{\nu_{x}\bigl(\Om_k(o) \cap \Om_l(x)\bigr)}
{\nu_{o}\bigl(\Om_k(o) \cap \Om_l(x)\bigr)} \quad \text{for all}\;\;
\xi \in \Om_k(o) \cap \Om_l(x)\,.
$$
\end{pro}

(Note that the nonempty ones among the sets $\Om_k(o) \cap \Om_l(x)$,
$k,l \ge 0$, form a partition of $\bd^*\T$ consisting of open-compact parts.)

\begin{proof}
Let $x, y \in \T$. By irreducibility, there is $n$ such
that $p^{(n)}(x,y) > 0$. Since
$$
\nu_{x} = \sum_{w} p^{(n)}(x,w) \,\nu_{w} 
\ge p^{(n)}(x,y) \,\nu_{y}\,,
$$
we find that $\nu_{y} \ll \nu_{x}$. 
Since $\nu_{x}$ is equidistributed and does not vanish on any $\Om_k(x)$,
the above formula for the Radon-Nikodym derivative is immediate. 
As the sets $\Om_k(o) \cap \Om_l(x)$ are open, it also follows that
$K(x,\cdot)$ is locally constant, whence continuous.
\end{proof}

\eqref{harm-meas} implies that for each $\xi \in \bd^*\T$, 
the function $x \mapsto K(x,\xi)$ is harmonic, and 
$K(o,\xi)=1$ by construction. 

We can compute $K(x,\xi)$ more explicitly in terms of the coefficients
$a_j$ of Lemma \ref{a_k},  taking into account Remark \ref{remark} when
the ``negative drift'' condition \eqref{hor-moment} holds.
The following is obtained by a lengthy, but completely straightforward discussion
of all possible relative positions of $\xi, x$ and $o$.
\begin{lem}\label{kernels}
Let $x \in \T$ and $\xi \in \bd^*\T$. Set $k= \up(o,\xi)$,
$l=\up(x,\xi)$ and $m = \hor(x)$, so that $k, l \ge 0$ and 
$m \in \Z$. Also, set $b(m) = 1$ in case A (non-negative drift) and
$b(m) = e^{c_0m}$ in case B (negative drift). Then
$$
K(x,\xi) = b(m)\, \frac{a_l}{a_k}\,q^{k-l+m} 
                 \left( \frac{q}{q-1}\right)^{\epsilon(k,l,m)}\,,
$$		 
where
$$
\begin{gathered}
\epsilon(k,l,m) = \sign(k-l+m)\sign(k)\sign(l)\,, \quad\text{except for}\\
\epsilon(0,l,m) = 1\,, \quad 1 \le l \le m\,\; (m > 0) \AND
\epsilon(k,0,m) = -1\,, \quad 1 \le k \le -m\,\; (m < 0)\,.
\end{gathered}
$$
In particular, the $P$-harmonic function $K(\cdot,\xi)$ is unbounded
for every $\xi \in \bd^*\T$.
\end{lem}
 
Our notation seems to indicate that each function $K(\cdot,\xi)$
is a Martin kernel. This is indeed true. 

\begin{thm}\label{sara} Suppose that $P$ is irreducible, and that
the respective moment conditions of cases A or B above hold.
Let $G(x,y)$ and $K(x,y)=G(x,y)/G(o,y)$
be the associated Green and Martin kernels. Then, with $K(x,\xi)$ 
as defined in Proposition \ref{derivatives}, we have
$$
\lim_{n \to \infty} K(x,y_n) = K(x,\xi)\,,
$$
whenever $(y_n)$ is a sequence of vertices that tends to 
$\xi \in \bd^*\T$ in the topology of the end compactification $\wh\T$.
\end{thm}

\begin{proof} We only need to consider case A. Indeed, in case B, 
when \eqref{hor-moment} holds, then $m(P^{\sharp}) < \infty$ 
and $\al(P^{\sharp} > 0)$. Since 
$K(x,y) = e^{c_0\hor(x)} K^{\sharp}(x,y)$, the result will
follow from case A.  

The proof is based on {\sc Brofferio} \cite{Bro}, Thm. 3.6.2--3. 
We briefly explain how that result has to be ``translated'' to our situation.
Let $\mu$ be as in \eqref{mu}, and let $\nu$ be the 
$\mu$-invariant measure on $\bd^*\T$ (unique up to normalization) according
to Proposition \ref{invariant}. 

Besides the group $\Ga=\AFF(\T)$ and the stabilizer $\Ga_{\! o}$, we also
need the \emph{horocyclic subgroup} $\Hor(\Ga) = \{ g \in \Ga : \hor(go)=0 \}$
of all elements of $\Ga$ that stabilize some (whence every) horocycle as
a set. Let $\KK=\KK_{\xi}$ denote the stabilizer of $\xi \in \bd^*\T$ in 
$\Hor(\Ga)$. It is a compact subgroup, and 
since $\Hor(\Ga)$ acts transitively on $\bd^*\T$, we can identify
$\bd^*\T$ with $\Hor(\Ga)/\KK$. Thus, we can lift $\nu$ to a
right-$\KK$-invariant measure $\bar \nu$ on $\Hor(\Ga)$.
More precisely, for each $\eta \in \bd^*\T$, let $g_{\eta}$ be
an element of $\Hor(\Ga)$ with $g_{\eta}\xi = \eta$. Also, let
$\la_{\KK}$ be the Haar measure of $\KK$, normalized with total mass $1$.
Then, for any Borel set $B \subset \Hor(\Ga)$,
$$
\bar\nu(B) = \int_{\bd^*\T} 
\la_{\KK}\bigl(g_{\eta}^{-1}\{g \in B:g\xi = \eta\}\bigr)\,
d\nu(\eta)\,.
$$
This is independent of the specific choice of $g_{\eta}$, but $\bar \nu$
does depend on $\xi$. In particualr, if the set $B$ is right-$\KK$-invariant,
then $\bar\nu(B) = \nu(B\xi)$.

Next, given $\xi \in \bd^*\T$, we can choose a ``shift'' $s = s_{\xi}$ in 
$\Ga$ that maps each element of the geodesic $\geo{\om\,\xi}$ to its 
successor on $\geo{\om\,\xi}$. Let $\la_{\langle s \rangle}$ be the counting
measure on the cyclic subgroup $\langle s \rangle$ of $\Ga$. 

We can consider both $\bar \nu$ and $\la_{\langle s \rangle}$ as measures
on the whole of $\Ga$, supported by the respective subgroups. We also
consider the potential associated with $\mu$, that is, the Radon measure
$\U = \sum_{n=0}^{\infty} \mu^{(n)}$, where $\mu^{(n)}$ is the $n$-th
convolution power of $\mu$. By transience, $\U$ is a Radon measure on $\Ga$.
Then Thm. 3.6 of \cite{Bro} says that in case A, with the proper choice 
of the normalizing constant $c(\mu)>0$,
\begin{equation}\label{U}
\U * \de_{g_n^{-1}} \to c(\mu)\cdot \bar\nu * \la_{\langle s \rangle}
\quad \text{vaguely, as $n \to \infty$,}
\end{equation}
whenever $(g_n)$ is a sequence in $\Ga$ such that $g_no \to \xi$
in the end topology. (As a matter of fact, \cite{Bro} states and
proves the ``inverse'' statement,  and we are applying
that result to the reflected measure $\check \mu$.)

We now translate this to our situation. Recall that the stabilizer
$\Ga_{\! o}$ is open and compact, and also recall the definition \eqref{mu}
of $\mu$ in terms of left Haar measure on $\Ga$. It implies that 
$\U(g\Ga_{\! o}) = \Delta(g^{-1})\U(\Ga_{\! o}g)$ for every $g \in \Ga$, where
$\Delta$ is the modular function of $\Ga$. It is known that for 
$g \in \Ga = \AFF(\T_q)$, the latter is $\Delta(g) = q^{\hor(go)}$, see e.g. 
\cite{CaKaWo}. For $y_n$ tending to $\xi$, we find $g_n \in \Ga$ such that 
$g_nx = y_n$. We obtain for the Green kernel of $P$
$$
G(o,y_n) = \U(g_n\Ga_{\! o}) 
= q^{-\hor(y_n)} \U*\de_{g_n^{-1}}(\Ga_{\! o})\,.
$$
Consequently, \eqref{U} implies
$$
q^{\hor(y_n)}G(o,y_n) \to 
c(\mu)\cdot \bar\nu * \la_{\langle s \rangle} (\Ga_{\! o}) = 
c(\mu)\cdot \bar\nu(\Ga_{\! o}) 
$$
In order to compute $\bar\nu(\Ga_{\! o})$, observe that $\Ga_{o}\!\KK$ is 
the disjoint union of $|\KK \!o|$ right $\Ga_{\! o}$-cosets, and that 
$|\KK\!  o| = |T_{k,k}|$, where $T_{k,k}$ is as in \eqref{Om_k} and
$k=\up(o,\xi)$. Since the measure $\bar\nu$ and the set
$\Ga_{\! o}\KK$ are right-$\KK$-invariant, 
$$
|T_{k,k}|\, \bar\nu(\Ga_{\! o}) = \bar\nu(\Ga_{\! o}\! \KK) 
= \nu(\Ga_{\! o}\! \KK\xi) = \nu(\Ga_{\! o}\xi) = \nu(\Om_k)\,.
$$
We obtain that $\bar\nu(\Ga_{\! o}) = \nu\bigl(\Om(\xi)\bigr)$,
as defined in \eqref{Om_k}.

Now we observe that for any $g \in \Ga$ and $y \in \T$, setting 
$\eta = g^{-1}\xi$, we have 
$$
\hor(g^{-1}y) = \hor(y) - \hor(go) \AND
\nu\bigl(\Om(\eta)\bigr) =  
q^{-\hor(go)} \de_g*\nu\bigl(\Om(\xi)\bigr)\,.
$$
Therefore, for any $g \in G$, setting $\eta = g^{-1}\xi\,$,
$$
\begin{aligned}
q^{\hor(y_n) - \hor(go)}G(go,y_n) 
&= q^{\hor(g^{-1}y_n)}G(o,g^{-1}y_n)\\
&\to c(\mu)\cdot \nu\bigl(\Om(\eta)\bigr) =  
c(\mu)\cdot q^{-\hor(go)} \de_g*\nu\bigl(\Om(\xi)\bigr)\,.
\end{aligned}
$$
If $x \in \T$ then we can find $g \in \Ga$ with $go=x$, whence
$$
q^{\hor(y_n)}G(x,y_n) \to 
c(\mu)\cdot \nu_{x}\bigl(\Om(\xi)\bigr)\,.
$$
From this, the result follows.
\end{proof}

\medskip

{\bf The minimal $P$-harmonic functions.} We can now determine the 
minimal harmonic functions. 
First of all we recall a well known fact that does not require any moment
condition; see e.g. \cite{Wbook}, Thm. 25.4.

\begin{pro}\label{minimal-Z}
Let the function $\varphi$ be defined as in \eqref{phi}. Then the minimal
$\wt P$-harmonic functions on $\Z$ are precisely the functions
$m \mapsto e^{cm}$, where $c \in \R$ is such that $\varphi(c)=1$.

Thus, if there is $c_0 \ne 0$ such that $\varphi(c_0) =1$ then the positive
$\wt P$-harmonic functions on $\Z$ are precisely the functions
$$
m \mapsto t + (1-t)e^{c_0m}\,,\quad t \in [0\,,\,1]\,.
$$
Otherwise, all positive $\wt P$-harmonic functions are constant.
\end{pro}

\begin{lem}\label{omega} If $(y_n)$ is a sequence in $\T$ tending to
$\om$ and such that $K(x,y_n) \to h(x)$ pointwise,
then $h$ depends only on $\hor(x)$, i.e., there is a $\wt P$-superharmonic
function $f$ on $\Z$ such that $h(x) = f\bigl(\hor(x)\bigr)$ for
all $x \in \T$.
\end{lem}

\begin{proof} By the hypothesis, $\up(o,y_n) \to \infty$. 
Thus, the proof is exactly as in \eqref{up-up}, see Figure 2 and the
subsequent lines. In general (unless $P$ has finite range), the limit
function is superharmonic, but not necessarily harmonic.
\end{proof}

\begin{thm}\label{main2}
Suppose that \emph{(i)} $\;m(P) < \infty$ and $\al(P) > 0\,$, or that 
\emph{(ii)} $\;m_{2+\ep}(P) < \infty$ and $\al(P) = 0\,$, or that
\emph{(iii)} condition \eqref{hor-moment} holds.

Then each function $K(\cdot,\xi)$, $\xi \in \bd^*\T$ is minimal 
$P$-harmonic.\\[6pt]\indent
In case \emph{(i)}, we have the following\\[4pt]
{\rm (a)} If the only solution to the equation $\varphi(c)=1$ is $c=0$, 
then the above are all minimal $P$-harmonic functions.\\[4pt]
{\rm (b)} Otherwise, if $\varphi(c_0) = 1$ for some $c_0 \ne 0$ (whence
$c_0 < 0$), the minimal $P$-harmonic functions are the above together 
with the function $h(x)=e^{c_0\hor(x)}$.\\[6pt]\indent
In cases \emph{(ii)} and \emph{(iii)}, the minimal $P$-harmonic functions are 
the above together with the constant function $h(\cdot) \equiv 1$.
\end{thm}

\begin{proof}
We start with a similar argument as in the proof of Theorem \ref{main1}.

First of all, note that there are positive harmonic functions that
are not constant on horocycles, namely the functions $K(\cdot,\xi)$
with $\xi \in \bd^*\T$. Therefore, there must be at least one minimal 
harmonic functions $h$ with the same porperty. Then there is a sequence
$(y_n)$ in $\T$ such that $K(\cdot,y_n) \to h$ pointwise.
By compactness of $\wh\T$, we may assume that $(y_n)$ converges in the
end topology to a point of $\xi \in \bd\T$. It cannot be $\xi=\om$,
because in that case $h$ would be constant on horocycles by Lemma \ref{omega}.
Thus, $\xi \in \bd^*\T$. Therefore, for this specific $\xi$, the function
$K(\cdot,\xi)$ is minimal harmonic.

Now let $\eta \in \bd^*\T$ be arbitrary. Then there is $g \in \Ga = \AFF(\T)$
such that $g\xi=\eta$. Note that we have the cocycle identity
\begin{equation}\label{cocycle}
K(x,g\xi) = K(g^{-1}x,\xi)/K(g^{-1}o,\xi)\,.
\end{equation}
Also, a function $h(x)$ is harmonic if and only if $h(gx)$ is harmonic.
Using these observations, it is a straightforward exercise that
$K(\cdot,\eta)$ is minimal as well. This proves the first part.

\smallskip
Suppose that $h$ is a minimal harmonic function distinct from all
$K(\cdot,\xi)$, $\xi \in \bd^*\T$. Then we must have, using the same
argument as at the beginning of the proof, that $h = \lim K(\cdot,y_n)$
for a sequence $(y_n)$ that tends to $\om$ in the end compactification.
Therefore $h$ is constant on horocycles, that is, we can write
$h(x) = f\bigl(\hor(x)\bigr)$, where $f$ is a $\wt P$-harmonic 
function on $\Z$. Minimality of $h$ with respect to $P$ implies 
$\wt P$-minimality of $f$. By Lemma \ref{minimal-Z}, 
$h(x) = e^{c\hor(x)}$, where $c$ satisfies $\varphi(c)=1$.  

\medskip

Now suppose we are in case (i). Then 
$$
\int K(\cdot,\xi)\,d\nu_{o}(\xi) = 1\,,
$$
whence the constant harmonic function $1$ is not minimal harmonic.
Thus, if $c=0$ is the only solution of $\varphi(c)=1$, then $h$ as above cannot
exist, and statement (a) holds.

Otherwise, we have to verify that $h(x)=e^{c_0\hor(x)}$ is indeed minimal.
As stated, we must have $c_0 < 0$, since the (at least one-sided)
derivative of $\varphi$ at $0$ is positive, and $\varphi$ is convex.

Suppose that $h$ is not minimal. Then the minimal harmonic functions are
precisely the $K(\cdot,\xi)$, $\xi \in \bd^*\T$. Thus, there is
a probability measure $\nu^h$ on $\bd^*\T$ such that 
$$
h(x) = \int_{\bd^*\T} K(x,\cdot) \,d\nu^h\,.
$$
A straightforward computation based on Lemma \ref{kernels} shows that 
for $r > k \ge 0$,
$$
\sum_{x \in T_{k,r}} K(x,\xi) = 
q^r\frac{a_0 + \dots + a_{r-1} + \kappa_k\,a_r}{a_k}
\quad\text{for all}\; \xi \in \Om_k\,,
$$
where $\kappa_k = 1$ if $k \ge 1$ and $\kappa_0=\frac{q-2}{q-1}$.
Therefore
$$
\begin{aligned}\left(\frac{q-1}{q}\right)^{\sign(k)}q^r\,e^{c_0(r-k)} 
&= \sum_{x \in T_{k,r}} f(x) 
\ge \int_{\Om_k} \left(\sum_{x \in T_{k,r}} K(x,\cdot)\right)\,d\nu^h\\
&\ge q^r \, \frac{a_0 + \dots + a_{r-1} + \kappa_k\,a_r}{a_k}\, \nu^h(\Om_k)\,.
\end{aligned}
$$
We see that 
$$
a_k \, e^{c_0(r-k)} \ge (a_0 + \dots + a_{r-1})\,\nu^h(\Om_k)
$$
for all $r > k$. Letting $r\to \infty$, the left hand side tends to $0$,
while the right hand side tends to $\nu^h(\Om_k)$. Thus, the measure
$\nu^h$ vanishes everywhere, a contradiction.

\medskip

In case (ii), the only solution of $\varphi(c)=1$ is $c=0$. The associated
harmonic function is $h(\cdot) \equiv 1$. It is known from \cite{CaKaWo} that 
in the case when $\al(P)=0$, the Poisson boundary is trivial, that is, all
bounded harmonic functions are constant. This amounts to minimality
of $h(\cdot) \equiv 1$. The same argument as used in case (i) shows that
there can be no further minimal $P$-harmonic functions.

\medskip 

Case (iii) is immediate by applying case (i.b) to $P^{\sharp}$,
since a function $h$ is minimal harmonic for $P$ if and only if
$h^{\sharp}(x) = e^{-c_0\hor(x)}h(x$ is minimal harmonic for
$P^{\sharp}$. Here, $c_0 > 0$ is the constant of \eqref{hor-moment}.
\end{proof}

\section{Conclusion}\label{conclusion}

We now return to the ``lamplighter setting'', where $P$ is a semi-isotropic,
irreducible transition matrix on $\DL(q,r)$. We write $P_1$ and $P_2$ for the 
projections of $P$ onto $\T^1=\T_q$ and  $\T^2 = \T_r$, respectively.
We can apply all results of the preceding \S \ref{tree} to each $P_i$. 
If $\wt\mu_i$ denotes the measure on $\Z$ that describes the projection 
$\wt P_i$ of $P_i$, i.e., $\wt p_i(k,l)=\wt\mu_i(l-k)$, then 
$\wt\mu_2(k) = \wt\mu_1(-k)$, whence $\varphi_2(c) = \varphi_1(-c)$ for the
associated functions according to \eqref{phi}. We shall stick to
$\varphi = \varphi_1$, which in terms of $P$ on $\DL$ is given by
$$
\varphi(c) = \sum_{x_1x_2 \in \DL} e^{c \hor(x_1)}\,p(o_1o_2,x_1x_2)\,,
$$
and 
$$
\al(P) = \al(P_1) = \sum_{x_1x_2 \in \DL} \hor(x_1)\,p(o_1o_2,x_1x_2)\,,
$$
if the latter series converges absolutely (so that $\al(P_2) = -\al(P_1)$).
The moments $m_t(P)$ are defined as in \eqref{moment}. In addition, we
also introduce the exponential moment
$$
m^{(c)}(P) = \sum_{x_1x_2 \in \DL} 
\Bigl(d(o_1,x_1)\,e^{c_+ \hor(x_1)} + d(o_2,x_2)\,e^{c_- \hor(x_2)}\Bigr)
p(o_1o_2,x_1x_2)\,,
$$
where $c_+ = \max\{c,0\}$ and $c_- = \min\{c,0\}$. The purpose of this
condition is the following. Suppose that
there is $c_0 \ne 0$ such that $\varphi(c_0)=1$ and $m^{(c_0)}(P) < \infty$.
If $c_0 > 0$ then $P_1$ satisfies  \eqref{hor-moment} on $\T^1$, and  
Theorem \ref{main2} applies to $P_1$. Also, $m_1(P_2) < \infty$ in that case,
whence $\al(P_2)$ exists, and it must be $\al(P_2) > 0$, since
$\al(P_1) = -\al(P_2)$ cannot be non-negative. Therefore, Theorem 
\ref{main2} also applies to $P_2$. If $c_0 < 0$, the situation is analogous,
with the roles of $P_1$ and $P_2$ exchanged.

In each case, we write $K_i(x_i,\xi_i)$ for the respective kernels 
on $\T^i$ according to Proposition \ref{derivatives} and Lemma \ref{kernels}.

\begin{thm}\label{main3} Let $P$ be an irreducible, semi-isotropic transition 
matrix on $\DL(q,r)$, and $P_i$ ($i=1,2$) its projections onto the trees
$\T^1=\T_q$ and $\T^2=\T_r$, respectively. Suppose that\/ {\rm (I)} 
$\;m_{2+\ep}(P) < \infty$ and $\al(P)=0\,$, or that\/ {\rm (II)} 
there is $c_0 \ne 0$ such that $\varphi(c_0)=1$ and $m^{(c_0)}(P) < \infty$.

\smallskip

Then each of the functions $x_1x_2 \mapsto K_i(x_i,\xi_i)$, where 
$\xi_i \in \bd^*\T^i$ ($i=1,2$) is a minimal $P$-harmonic function on
$\DL(q,r)$.

\smallskip

In case (I), the minimal harmonic functions are the above together with the 
constant function $h(\cdot) \equiv 1$.  In case (II), the above are all
minimal harmonic functions.
\end{thm}

\begin{proof} Combining Theorem \ref{main1}(a) with Theorem \ref{main2},
we conclude that each minimal $P$-harmonic function must be of the form 
$x_1x_2 \mapsto K_i(x_i,\xi_i)$ with $\xi_i \in \bd^*\T^i$ ($i=1,2$),
or $x_1x_2 \mapsto e^{c\hor(x_1)}$ with $\varphi(c)=1$. 

Since there are $P$-harmonic functions that depend only on $x_i$, but 
are not of the form $x_1x_2 \mapsto f\bigl(\hor(x_i)\bigr)$,
there also must me a minimal harmonic function with these properties.
By the above, it must be a kernel $x_1x_2 \mapsto K_i(x_i,\xi_i)$.
Hence, there is at least one $\xi_i \in \bd^*\T^i$ such that
$x_1x_2 \mapsto K_i(x_i,\xi_i)$ is minimal $P$-harmonic. Using the
same cocycle argument as below \eqref{cocycle} in the proof of Theorem \ref{main2},
we obtain that \emph{all} $\xi_i \in \bd^*T^i$ ($i=1,2$) give rise
to a minimal $P$-harmonic function on $\DL$.

In case (I), by Theorems \ref{main1} and \ref{main2} the only other 
candidate for being a minimal $P$-harmonic function is the constant 
function $h(\cdot) \equiv 1$. The latter is indeed minimal:
by Theorem \ref{main1}, every postive 
bounded $P$-harmonic function is of the form $h(x_1x_2) = h_1(x_1)+h_2(x_2)$, 
where each $h_i$ must be bounded $P_i$-harmonic, whence constant by
Theorem \ref{main2}, as $\al(P_i)=0$. Now recall that the constant function
$1$ is minimal if and only if all bounded harmonic functions
are constant.

In case (II), the only candidates besides the kernels $K_i(\cdot,\xi_i)$
for being minimal $P$-harmonic are the functions $x_1x_1 \mapsto e^{c\hor(x_1)}$
with $c=0$ and $c=c_0$. Regarding $c=0$, we know that the constant function
$h(\cdot) \equiv 1$ is not minimal, since it is not minimal for $P_i$ on $\T^i$,
where $i$ is the index for which $\al(P_i) > 0$. Analogously, if we define
$P^{\sharp}$ by 
$p^{\sharp}(x_1x_2,y_1y_2) = p(x_1x_2,y_1y_2)\,e^{c_0(\hor(y_1)-\hor(x_1))}$,
then $h (\cdot) \equiv 1$ is not minimal for $P^{\sharp}$, whence
$x_1x_2 \mapsto e^{c_0\hor(x_1)}$ is not minimal for $P$.
\end{proof}

\end{document}